\documentclass{article}
\usepackage{emsnewsletter}
\usepackage{youngtab}
\usepackage{epsfig}
\usepackage{color}
 
\newtheoremstyle{exercise}
  {3pt}
  {6pt}
  {}
  {}
  {\bfseries}
  {:}
  { }
   {}
\theoremstyle{exercise}
\newtheorem{exercise}[theorem]{Exercise}
\newtheoremstyle{exercises}
  {3pt}
  {6pt}
  {}
  {}
  {\bfseries}
  {:}
  {\newline}
   {}

\theoremstyle{exercise}
\newtheorem{exercises}[theorem]{Exercises}

\input epsf
\def\boxit#1{\vbox{\hrule height1pt\hbox{\vrule width1pt\kern3pt
  \vbox{\kern3pt#1\kern3pt}\kern3pt\vrule width1pt}\hrule height1pt}}

\def\trank{\text{rank}}

\def\bv{\bold v}

\def\BC{\mathbb C}
\def\BR{\mathbb R}
\def\BP{\mathbb P}
\def\pp#1{\mathbb P^{#1}}

\def\pp#1{{\mathbb P}^{#1}}
\def\tdim{{\rm dim}}

\def\hd{,...,}
\def\ww{\wedge}
\def\upperp{{}^\perp}

\def\11{\mathbf 1}

\def\l{\lambda}
\def\a{\alpha}

\def\b{\beta}
\def\g{\gamma}
\def\s{\sigma}

\def\d{\delta}

\def\ot{{\mathord{ \otimes } }}
\def\op{{\mathord{\,\oplus }\,}}
\def\otc{{\mathord{\otimes\cdots\otimes}\;}}

\def\ra{{\mathord{\;\rightarrow\;}}}

\def\dim{{\rm dim}\;}
\def\La#1{\Lambda^{#1}}

\def\op{\oplus}
\def\BZ{\Bbb Z}

\def\op{\oplus}

\def\s{\sigma}

\def\a{\alpha}
\def\b{\beta}

\def\g{\gamma}
\def\l{\lambda}

\def\FS{\mathfrak  S}

\def\ol{\overline}

\def\BP{\mathbb  P}
\def\BC{\mathbb  C}

\def\pp#1{\mathbb  P^{#1}}

\def\BR{\mathbb  R}

\def\hd{, \hdots ,}

\def\La#1{\Lambda^{#1}}

\def\pp#1{\mathbb  P^{#1}}

\def\ra{\rightarrow}

\def\tdet{\operatorname{det}}

\def\tperm{\operatorname{perm}}

\def\tend{\operatorname{End}}

\def\tdim{\operatorname{dim}}

\def\tmin{\operatorname{min}}

\def\thom{\operatorname{Hom}}
\def\trank{\operatorname{rank}}

\def\upperp{{}^{\perp}}

\def\ww{\wedge}

\def\be{\begin{equation}}
\def\ene{\end{equation}}

\def\tsgn{{\rm{sgn}}}

\DeclareMathOperator{\tlog}{log}

\def\Mn{M_{\langle  n \rangle}}

\def\dual{{^\vee}}

\def\vp{{\bold V\bold P}}
\def\vnp{{\bold V\bold N\bold P}}\def\p{{\bold P}}
\def\np{{\bold N\bold P}}

\def\G{\Gamma}

\def\tzeros{{\rm Zeros}}

\def\Mn{M_{\langle \nnn \rangle}}\def\Mone{M_{\langle 1\rangle}}

\def\trank{{\mathrm {rank}}}

\def\tmult{{\rm mult}}

\def\Dual{{\mathcal D}}

\def\nnn{\bold n}

\def\bv{\bold v}

\def\Det{{\mathcal Det}}\def\Perm{{\mathcal Perm}}
\def\hDet{ {  \mathcal {\hat D}et} }\def\hPerm{{\mathcal {\hat P}  erm}}
\def\hzeros{ {   {\hat Z}eros} }

\author{J.M. Landsberg  (Texas A\&M University,  College Station, TX,  USA)}
\title{An introduction to geometric complexity theory}

\begin{document}
   
\maketitle

\begin{abstract}
I survey   methods from differential geometry,
 algebraic geometry and representation theory relevant for the permanent v. determinant
 problem from computer science, an algebraic analog of the $\p$ v. $\np$ problem.
\end{abstract}

\section{Introduction}\label{sec:intro}The purpose of this article is to introduce mathematicians to   uses  
of geometry in complexity theory.   I  focus on a central question:  
 the {\it Geometric Complexity Theory}  version of L. Valiant's conjecture comparing the complexity of the
 permanent and determinant polynomials, which is an algebraic
 variant of the  $\p \neq \np$ conjecture. Other problems in complexity such as {\it matrix rigidity}
 (see \cite{MR2870721,GHILrigid,aluffideg}) and the {\it complexity of matrix multiplication} (see, e.g.,  \cite{MR2383305}) have been treated with similar geometric  methods.
 

\section{History}\label{historysect}
\subsection{1950's Soviet Union}
A traveling saleswoman   needs to
visit 20 cities;  Moscow, Leningrad, Stalingrad,... Is there a route that
can be taken traveling less than 10,000km?

 Essentially the only known method to determine the answer is a brute force search through all possible
 paths.  
 The number of paths to check grows exponentially with the number of cities to visit.
  Researchers in the Soviet Union asked:
 {\it Is this brute force search avoidable?}  I.e.,  are there any algorithms that are
 significantly better than the na\"\i ve one?

 A possible  cause for hope is that if someone proposes a route, it is very easy to check if
 it is less than 10,000km (even pre-Google).

\subsection{1950's Princeton NJ}
In a letter to von Neumann (see \cite[Appendix]{Sipser}) G\"odel
  attempted  to quantify   the apparent difference between intuition
    and systematic problem solving.
For example, is it really significantly easier to verify a  proof than to
write one?

\subsection  {1970's: Precise versions of these questions}
These ideas evolved to a precise conjecture posed by  
 Cook  (preceded by work of Cobham, Edmonds,  Levin,  Rabin, Yablonski, and the above-mentioned question of G\"odel):
 
 Let $\p$ denote the class of problems that are \lq\lq easy\rq\rq\  to solve.\footnote{Can be solved 
 on a Turing machine in time polynomial with respect to the size of the input data.}
 
 Let $\np$ denote the class of problems that are \lq\lq easy\rq\rq\  to
 verify (like the traveling saleswoman problem).\footnote{A proposed solution can be verified in polynomial time.}
 
 \begin{conjecture}\cite{cook1971complexity,karp1972reducibility} $\p\neq\np$.\end{conjecture}

 \subsection  {Late 1970's: L. Valiant, algebraic variant}
 A {\it bipartite graph} is  a graph with two sets of vertices and edges   joining vertices from
one set to the other. A {\it perfect matching} is a subset of the edges such that each vertex shares
an edge from the subset
  with exactly one other vertex.
 
A  standard  problem in graph theory, for which the only known algorithms are exponential in the size of
the graph, is  to count the number of perfect matchings of a bipartite graph.

\begin{figure}[!htb]\begin{center}\label{biparga}
\includegraphics[scale=.3]{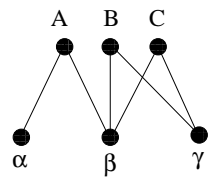}
\caption{\small{A bipartite graph, Vertex sets are $\{A,B,C\}$ and $\{\a,\b,\g\}$. }}  
\end{center}
\end{figure}

\begin{figure}[!htb]\begin{center}
\includegraphics[scale=.3]{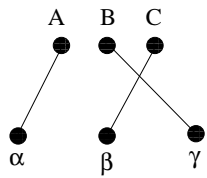}\ \ \ \ \ \ \ \ \includegraphics[scale=.3]{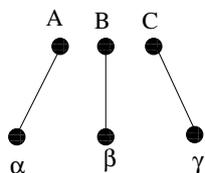}
\caption{\small{Two perfect matchings of the graph from Figure 1. }}  
\end{center}
\end{figure}

This count
  can be computed by evaluating a polynomial as follows:
  To a bipartite graph $\G$ one associates an incidence matrix $X_{\G}=(x^i_j)$, where
$x^i_j=1$ if an edge joins the vertex $i$ above to the vertex $j$ below and is
zero otherwise. For example the graph of Fig. 1 has incidence matrix
$$
X_{\G}=
\begin{pmatrix} 1&1&0\\ 0&1&1\\0&1&1\end{pmatrix}.
$$

A perfect matching corresponds to a set of entries
$\{ x^{1}_{j_1}\hd x^n_{j_n}\}$ 
with all $x^i_{j_i}=1$ 
and
$(j_1\hd j_n)$ is a permutation of
$(1\hd n)$. Let $\FS_n$ denote the group of permutations
of the elements $(1\hd n)$.

Define the {\it permanent} of an $n\times n$
matrix $X=(x^i_j)$ by
\be 
perm_n(X):=\sum_{\s\in \FS_n}  x^1_{\s(1)}x^2_{\s(2)}
\cdots x^n_{\s(n)}.
\ene 
Then $\tperm(X_{\G})$ equals the number of perfect matchings of $\G$.

For example, 
$
\tperm_3 \begin{pmatrix} 1&1&0\\ 0&1&1\\0&1&1\end{pmatrix} =2.
$

 A fast algorithm to compute the permanent  would give a fast algorithm to
 count  the number of perfect matchings of a bipartite graph.

\smallskip
While   it may not be  easy to evaluate,  the polynomial 
  $\tperm_n$ is relatively easy to write down  compared with
 a random polynomial of degree $n$ in $n^2$ variables in the following sense:

Let $\vnp$ be the set of sequences of polynomials that are \lq\lq easy\rq\rq\ to write 
down.\footnote{Such sequences are obtained from sequences in $\vp$ (defined in the
following paragraph) by
\lq\lq projection\rq\rq\ or \lq\lq integration over the fiber\rq\rq\ where
one averages the polynomial over a subset of its variables specialized to $0$ and $1$.}  Valiant showed \cite{MR526203} that the permanent
is {\it complete} for the class $\vnp$, in the sense that $\vnp$ is the class of  all polynomial sequences $(p_m)$, where $p_m$ has degree $m$ and involves a number of variables
polynomial in $m$, 
such that  there is a polynomial $n(m)$ and  $p_m$ 
is an {\it affine linear projection} of   $\tperm_{n(m)}$
  as defined   below.
Many problems from graph theory, 
combinatorics, and  statistical physics (partition functions)  are in $\vnp$. A good way to think of $\vnp$ is as the
class of sequences of polynomials that can be written down explicitly.\footnote{Here one must take a narrow view of
explicit- e.g. restrict to integer coefficients that are not \lq\lq too large\rq\rq .}

Let $\vp$  be the set of sequences of polynomials that are \lq\lq easy\rq\rq\ to 
compute.\footnote{Admit a polynomial size arithmetic circuit, and polynomially bounded degree see e.g.  \cite[\S 21.1]{BCS}}  For example, one can compute the determinant of an $n\times n$ matrix
quickly, e.g., using Gaussian elimination, so the sequence $(\tdet_n)\in \vp$. 
 Most problems from linear algebra (e.g., inverting a matrix, computing its determinant, multiplying
matrices) are in $\vp$.

The standard formula for the easy to compute determinant polynomial is
\be 
det_n(X):=\sum_{\s\in \FS_n} \tsgn(\s) x^1_{\s(1)}x^2_{\s(2)}
\cdots x^n_{\s(n)}.
\ene 

Here $\tsgn(\s)$ denotes the sign of the permutation $\s$.
 
 Note that 
\begin{align*}
\tperm_2\begin{pmatrix} y^1_1&y^1_2\\ y^2_1 & y^2_2\end{pmatrix}&= y^1_1y^2_2 + y^1_2y^2_1 
\\
&= \tdet_2\begin{pmatrix} y^1_1&-y^1_2\\ y^2_1 & y^2_2\end{pmatrix}.
\end{align*}

On the other hand, Marcus and Minc \cite{MR0147488},  building on work of  P\'olya and Szeg\"o (see  \cite{MR922386}),  
proved that one
could not  express $\tperm_m(Y)$ 
as a size $m$ determinant of a matrix whose entries are
affine linear functions of the variables $y^i_j$ when $m>2$.  
This raised the question that   perhaps   the permanent of an $m\times m$ matrix 
could be expressed as a slightly larger determinant.
More precisely, we say $p(y^1\hd y^M)$ is an {\it affine linear projection}
\index{affine linear projection}  of $q(x^1\hd x^N)$, if
there exist affine linear functions $x^{\a}(Y)=x^\a(y^1\hd y^M)$ such that
$p(Y)=q(X(Y))$. 
For example, B. Grenet  \cite{Gre11} observed that
\be\label{grenet3}
\tperm_3(Y)=\tdet_7 \begin{pmatrix} 0& y^1_1&y^2_1&y^3_1&0&0&0\\
0&1&0&0&y^3_3&y^2_3&0\\
0&0&1&0&0&y^1_3&y^3_3\\
0&0&0&1&y^1_3&0&y^2_3\\
y^2_2&0&0&0&1&0&0\\
y^3_2&0&0&0&0&1&0\\
y^1_2&0&0&0&0&0&1
\end{pmatrix}.
\ene  
Recently \cite{2015arXiv150502205A} it was shown  that $\tperm_3$ cannot be realized as an affine linear projection of $\tdet_n$ for $n\leq 6$,
so \eqref{grenet3} is optimal.

Valiant   showed that if $n(m)$ grows exponentially with respect to $m$, then there exist affine linear functions
$x^i_j(y^s_t)$ such that
$\tdet_n(X(Y))=\tperm_m(Y)$. (Grenet strengthened this   to
show explicit expressions when $n= 2^m-1$ \cite{Gre11}. See \cite{LRpermdet} for a discussion of the geometry of these algorithms and a proof of their optimality
among algorithms with symmetry.)
Valiant also conjectured that one cannot do too much better:

\begin{conjecture} [Valiant \cite{MR526203}]\label{valperdet}  Let $n(m)$ be a polynomial of $m$. Then there exists an $m_0$ such that
for all $m>m_0$, there do not exist affine linear functions $x^i_j(y^s_t)$ such
that $\tperm_m(Y)=\tdet_n(X(Y))$.
\end{conjecture}

\begin{remark} The original $\p\neq\np$ is viewed as  completely out
of reach.  Conjecture \ref{valperdet}, which would be  implied
by    $\p\neq\np$ is viewed as a more tractable substitute. 
\end{remark}

To keep track of progress on the conjecture, for a polynomial
$p=p(Y)$, let $dc(p)$ denote the
smallest $n$ such that there exists 
an affine linear map $X(Y)$ satisfying $p(Y)=\tdet_n(X(Y))$. Then Conjecture \ref{valperdet} says that
$dc(\tperm_m)$ grows faster than any polynomial. Since the conjecture is expected to be
quite difficult, one could try to prove any lower bound  on $dc(\tperm_m)$.  Several linear bounds on $dc(\tperm_m)$ were shown
\cite{MR0147488,MR910987,MR1032157} 
 with the current world record  the quadratic bound $dc(\tperm_m)\geq \frac{m^2}2$
 \cite{MR2126826}. (Over finite fields one has the same bound by \cite{MR1032157}. Over $\BR$, one has $dc_{\BR}(\tperm_m)\geq m^2-2m+2$ \cite{DBLP:journals/corr/Yabe15}.)
The  state of the art was  obtained with local differential geometry, as described in 
\S\ref{valconjdg}.
 
\begin{remark} There is nothing special about the permanent for this conjecture: it would be sufficient to
show any explicit  (in the sense of $\vnp$  mentioned above) sequence of polynomials $p_m$ has $dc(p_m)$ growing faster
than any polynomial. The dimension of the set of affine linear projections
of $\tdet_n$ is roughly $n^4$, but   the dimension of the space of homogeneous polynomials of degree
$m$ in $m^2$ variables grows almost like $m^m$, so a   random    sequence will
have exponential $dc(p_m)$. Problems in computer science  to find an explicit object satisfying a property that  a random
one satisfies are called  {\it trying to find hay in a haystack}. 
\end{remark}
 
\subsection{Coordinate free version}
To facilitate the use of geometry, we   get rid of coordinates.    Let $\tend(\BC^{n^2})$ denote the space of   linear maps
$\BC^{n^2}\ra \BC^{n^2}$, which acts on the space of homogeneous
polynomials of degree $n$ on $\BC^{n^2}$, denoted $S^n\BC^{n^2 *}$ (where the $*$ is
used to indicate the dual vector space to $\BC^{n^2}$),  as follows:
for $g\in \tend(\BC^{n^2})$ and $P\in S^n\BC^{n^2 *}$, 
the polynomial
$ g\cdot P$ is defined by 
\be\label{gactonp}(g\cdot P)(x):=P(g^T x).
\ene
Here $g^T$ denotes the transpose of $g$.   (One takes the transpose matrix in order that
$g_1\cdot(g_2 \cdot P)=(g_1g_2)\cdot P$.)

In \cite{MS1} they introduced
  {\it padding}; adding a homogenizing variable so all objects live in the same
ambient space, in order to  deal with   linear functions instead of
affine linear functions. Let  
$\ell$ be a new variable, so $\ell^{n-m}\tperm_m(y)\in S^n\BC^{(m^2+1)*}$. Then   $\tperm_m(y)$ is expressible as an $n\times n$
determinant whose entries are affine linear combinations of the  $y^i_j$   if and only if $\ell^{n-m}\tperm_m$ is expressible
as an $n\times n$ determinant whose entries are {\it linear} combinations of  the variables $y^i_j,\ell$.

 Consider any linear inclusion
$  \BC^{(m^2+1) *}\ra \BC^{n^2  *}$, so in particular $\ell^{n-m}\tperm_m\in S^n\BC^{n^2  *}$.
Then  
\be\label{valreprhase}
dc(\tperm_m)\leq n \Leftrightarrow
\ell^{n-m}\tperm_m \in \tend(\BC^{n^2})\cdot \tdet_n.
\ene

Conjecture \ref{valperdet} in this language is:

\begin{conjecture} [Valiant \cite{MR526203}]   Let $n(m)$ be a polynomial of $m$. Then there exists an $m_0$ such that
for all $m>m_0$,  $\ell^{n-m}\tperm_m \not\in \tend(\BC^{n^2})\cdot \tdet_n$, equivalently
$$\tend(\BC^{n^2})\cdot\ell^{n-m}\tperm_m \not\subset  \tend(\BC^{n^2})\cdot \tdet_n.
$$
\end{conjecture}

\section{Differential geometry and the state of the art regarding Conjecture \ref{valperdet}}
\label{valconjdg}
The best result pertaining to  Conjecture \ref{valperdet} comes from local  differential geometry: the study of Gauss maps.
 
  \subsection  {Gauss maps}
Given a surface in $3$-space,  form its {\it Gauss map}  
by mapping a point of the surface   to its unit normal vector on the unit sphere as in Figure 3.

\begin{figure}[!htb]\begin{center}\label{coolgaussmap}
\includegraphics[width=9.6cm]{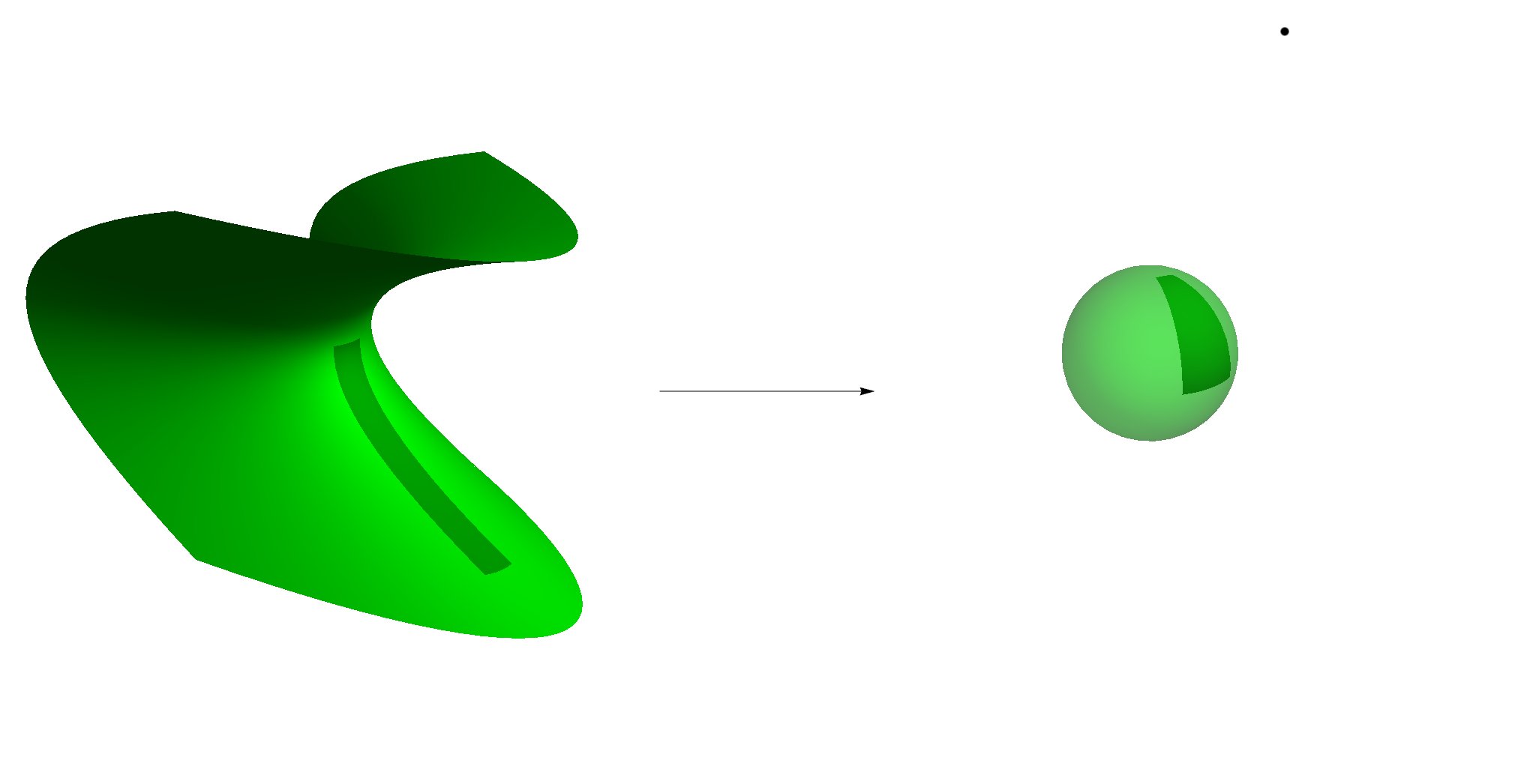} 
\caption{\small{The shaded area of the surface maps to the shaded area of the sphere.  }}  
\end{center}
\end{figure}

A normal vector to a surface $X$ at $x$ is one perpendicular to the tangent space $T_xX
\subset \BR^3$.
This Gauss image can be defined without  the use of an inner product if one instead takes the
union of all {\it conormal lines}, where a conormal vector to $X\subset \BR^3$ is one in the dual
space $\BR^{3*}$ that annhilates the tangent space $T_xX$. One loses   qualitiative  information, however
one still has the information of the {\it dimension} of the Gauss image.

This dimension will drop if through all points of the surface there is a curve 
along which the
tangent plane is constant. For example, 
   if $M$ is a   cylinder, i.e., the union of
lines in three space perpendicular to a plane curve, the Gauss image is a curve:

\begin{figure}[!htb]\begin{center} 
\includegraphics[height=3.7cm]{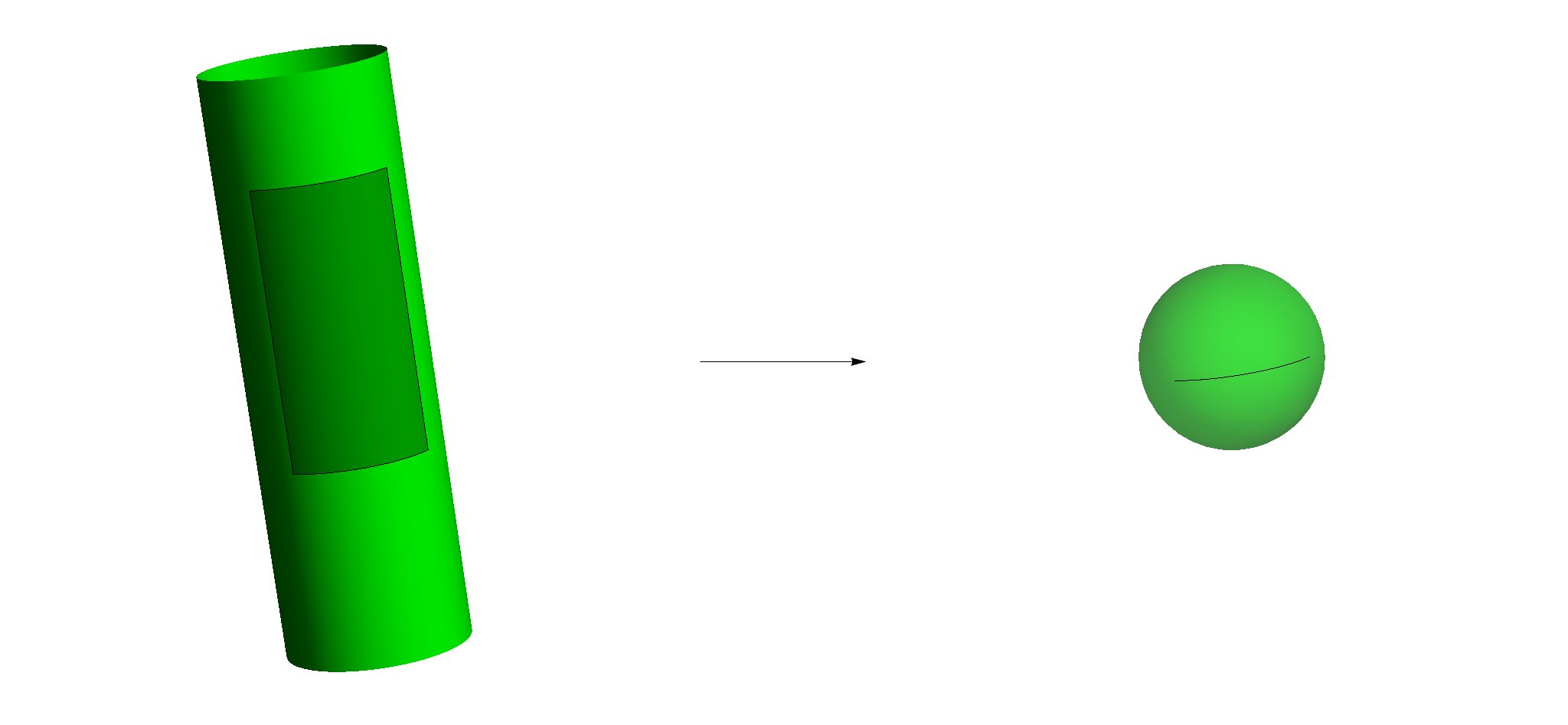} 
\caption{\small{Lines on the cylinder are collapsed to a point.}}  
\end{center}
\end{figure}

The extreme case is when  the surface is a plane, then its Gauss image is just a point.

A classical theorem in the geometry of surfaces in three-space classifies 
surfaces with degenerate Gauss image.   I state it in the algebraic category for what comes next  (for   $C^{\infty}$ versions see, e.g., \cite[vol. III, chap. 5]{MR532832}). 
One may view projective space $\pp 3$ as affine space with a plane added at infinity. From this perspective a cylinder is
a cone with vertex at infinity.

\begin{theorem} [C. Segre \cite{csegre}]
If $X^2\subset \pp 3$ is an algebraic surface
whose Gauss image is not two-dimensional,  then $X$ is one of:
\begin{itemize}
\item The   closure of the union of
points on  tangent lines to a space curve.
\item A generalized cone, i.e., the points on the union of lines connecting a fixed point to a plane curve. 
\end{itemize}
\end{theorem}

\includegraphics[width=4.3cm, angle=60]{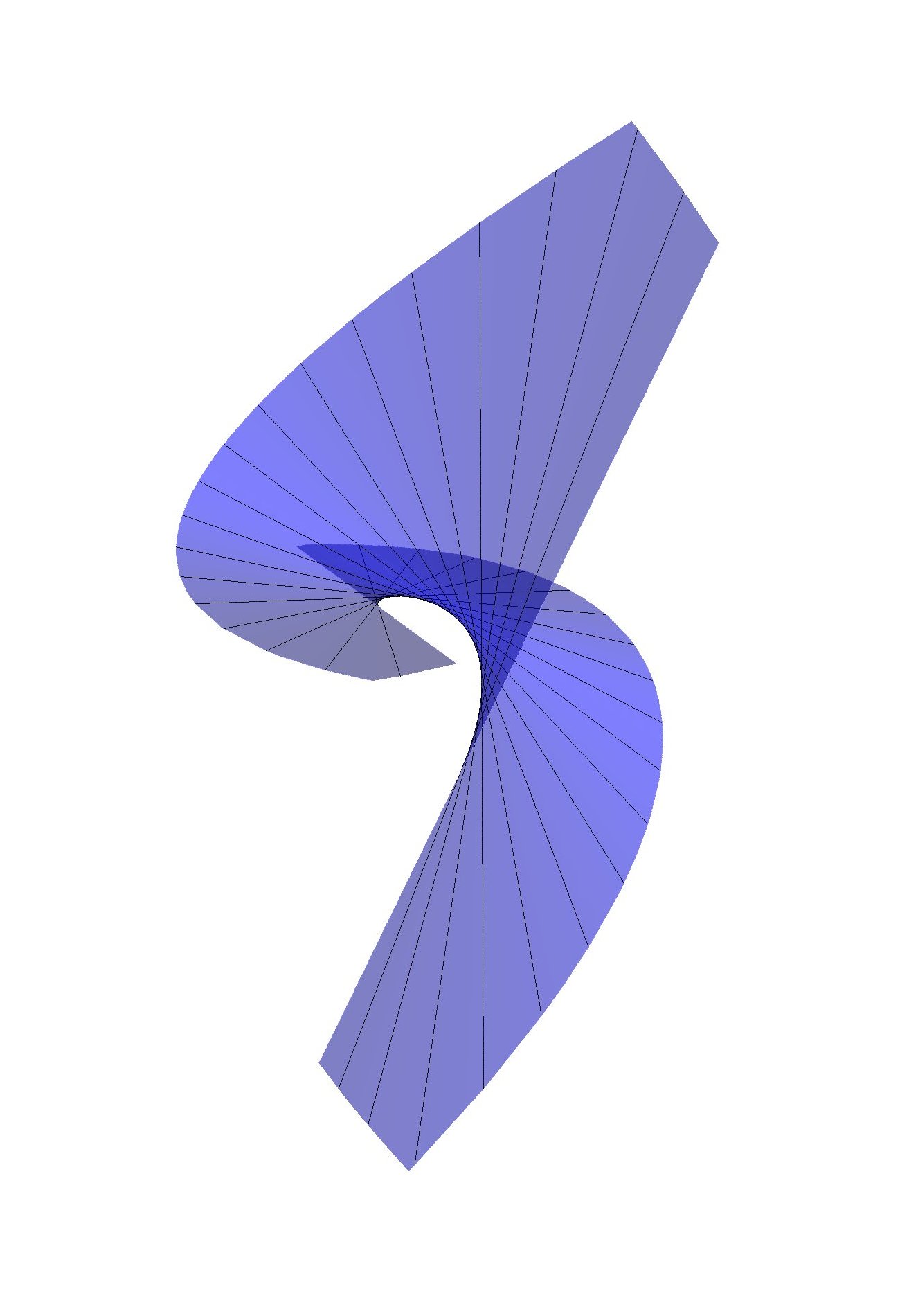}

\includegraphics[height=4.1cm]{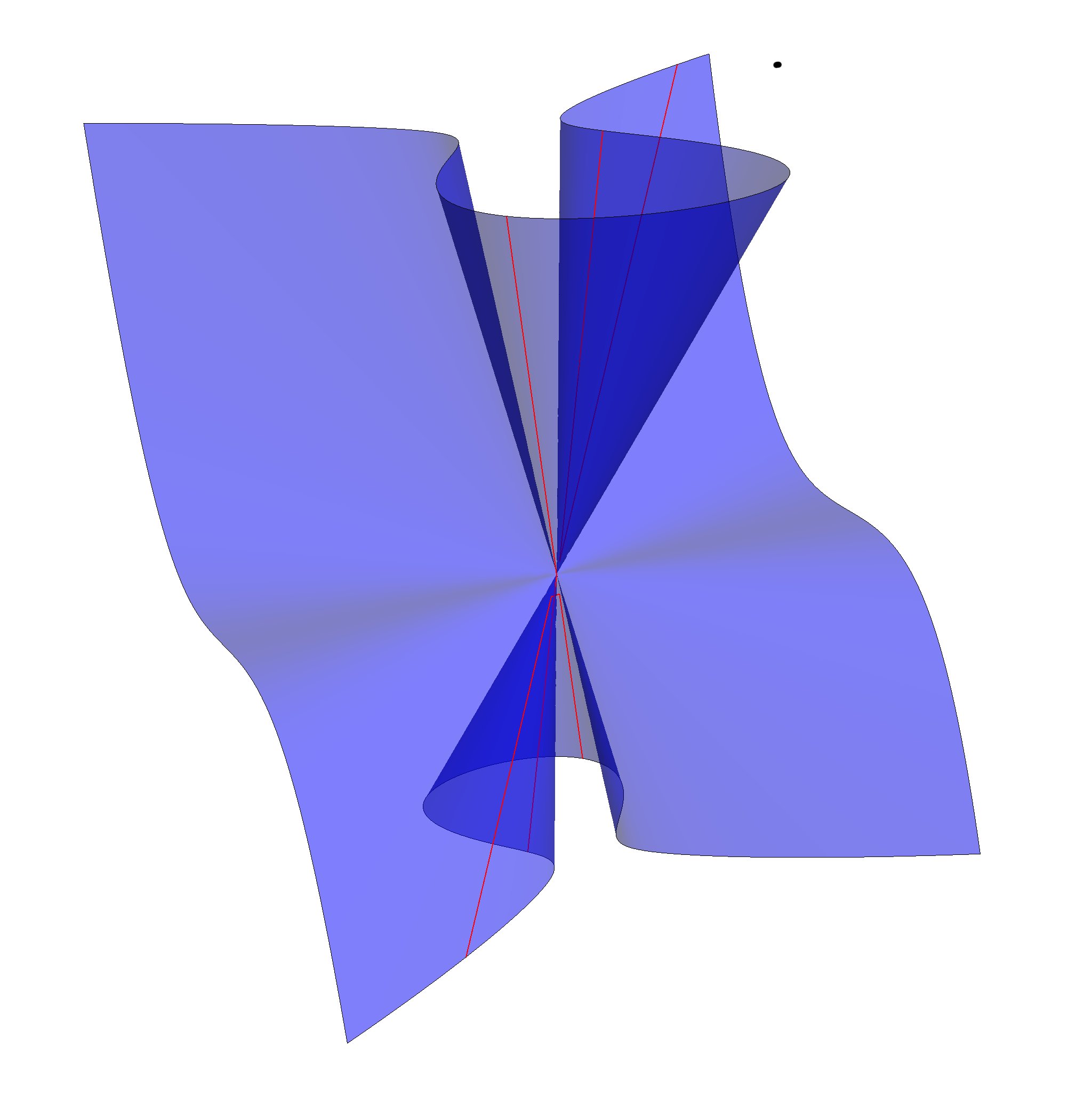}

Notice that in the first picture, the tangent plane along a ray of the curve is constant, and in the second
case the tangent plane is constant along the lines through  the vertex.

One can extend the notion of Gauss map to hypersurfaces of arbitrary dimension,
and to hypersurfaces defined over the complex numbers. 
The union of tangent rays to a curve generalizes to the case of {\it osculating varieties}.
One can also take cones with vertices larger than a point.

  \subsection  {What does this have to do with complexity theory?}
The hypersurface $\{\tdet_n(X)=0\}\subset \BC^{n^2}$ has
a very degenerate Gauss map. To see this,  consider the matrix
$$z=
\begin{pmatrix}
1  & & & \\ & \ddots & & \\ & & 1& \\ & & & 0\end{pmatrix} \in \{ \tdet_n=0\}.
$$
The  tangent space to $\{ \tdet_n=0\}$ at $z$, and the conormal space
(in the dual space of matrices) are respectively    
$$T_z\{ \tdet_n=0\}= \begin{pmatrix}
*  & *& *& *\\ \vdots  &  \ddots & \vdots  & \vdots  \\ * & * & *  & * \\ *& * & * & 0\end{pmatrix},
\ \ \ \ \ \ \
N^*_z\{ \tdet_n=0\}= \begin{pmatrix}
0  & 0& 0& 0\\ \vdots  &  \ddots & \vdots  & \vdots  \\ 0 & 0 & 0  & 0 \\ 0& 0 &  0 & *\end{pmatrix}.
$$
But any rank $n-1$ matrix whose non-zero entries all lie in the upper left $(n-1)\times (n-1)$ submatrix
will have the same tangent space!   
Since any smooth point of $\{\tdet_n=0\}$ can be moved to  $z$ by a change of
basis, we conclude that the tangent hyperplanes to  $\{\tdet_n=0\}$ are parametrized by the rank one matrices,
the space of which has dimension $2n-1$ (or $2n-2$ in projective space), because they
are obtained by multiplying a column vector by a row vector. In fact,
$\{\tdet_n=0\}$ may be thought of as an osculating variety of the variety of rank one matrices (e.g., the union of
tangent lines to the union of tangent lines... to the variety of rank one matrices). 

On the other hand, a direct calculation shows that the permanent hypersurface $\{\tperm_m=0\}\subset \pp{m^2-1}$
has a non-degenerate Gauss map (see \S\ref{lowerdcp}), so when one 
includes $\BC^{m^2}\subset \BC^{n^2}$,
   the   equation $\{\tperm_m=0\}$ becomes an equation in a space of $n^2$ variables
that only uses $m^2$ of the variables, one gets
a cone with vertex $\pp{n^2-m^2-1}$ corresponding to the unused variables, in particular, the Gauss image will   have dimension $m^2-2$.

If one makes an affine linear substitution $X=X(Y)$, the Gauss map of 
$\{\tdet(X(Y))=0\}$
will be at least as degenerate as the Gauss map of $\{\tdet(X)=0\}$. Using this, one obtains:

\begin{theorem}[Mignon-Ressayre \cite{MR2126826}] \label{mrthm} If $n(m)<\frac {m^2}2$, then there do not exist affine linear functions $x^i_j(y^s_t)$ such
that \linebreak
$\tperm_m(Y)=\tdet_{n}(X(Y))$. I.e., $dc(\tperm_m)\geq \frac {m^2}{2}$. 
\end{theorem}

\section{Algebraic geometry and Valiant's conjecture} \label{agandvalsect}
A possible path to show $\tperm_m(Y)\neq \tdet_n(X(Y))$   is to look for    a  polynomial  whose
zero set contains all polynomials of the form $\tdet_n(X(Y))$, and show that
$\tperm_m$ is not in the zero set.

\subsection{Polynomials} \label{LMRvague}Algebraic geometry is the study of zero sets of polynomials. 
In our situation, we need polynomials on spaces of polynomials.
More precisely, if 
$$P(x_1\hd x_N)=\sum_{1\leq i_1\leq \cdots\leq  i_d\leq N} c_{i_1\hd i_d}x_{i_1}\cdots x_{i_d}
$$
is a homogeneous polynomial of degree $d$ in $N$ variables, we work with  polynomials in the
coefficients $c_{i_1\hd i_d}$, where   these coefficients provide  coordinates
on the vector space $S^d\BC^N$ of all homogeneous polynomials of degree $d$ in $N$ variables.

The starting point of {\it  Geometric Complexity Theory}
is the plan  to prove Valiant's conjecture by finding a sequence of polynomials $P_m$ vanishing on all affine-linear projections
  of $\tdet_{n(m)}$ when $n$ is a polynomial in $m$ such that 
  $P_m$ does not vanish on $\tperm_m$.

\subsection{Disadvantage of algebraic geometry?} 
  The zero set of all polynomials vanishing on 
$$S:=\{ (z,w) \mid w=0,\ z\neq 0\}\subset \BC^2,$$ 

\begin{figure}[!htb]\begin{center}
\includegraphics[scale=.3]{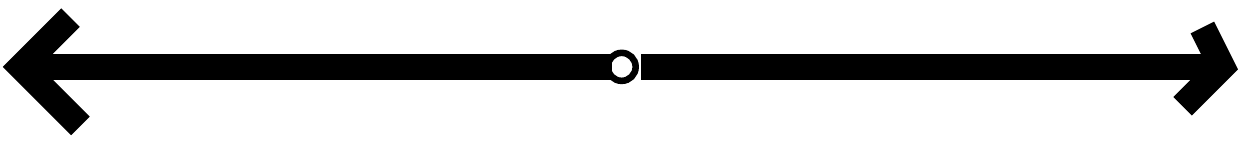}
\end{center}
\end{figure}
is the line  $$\{(z,w) \mid  w=0\}\subset \BC^2.$$  

\begin{figure}[!htb]\begin{center}
\includegraphics[scale=.3]{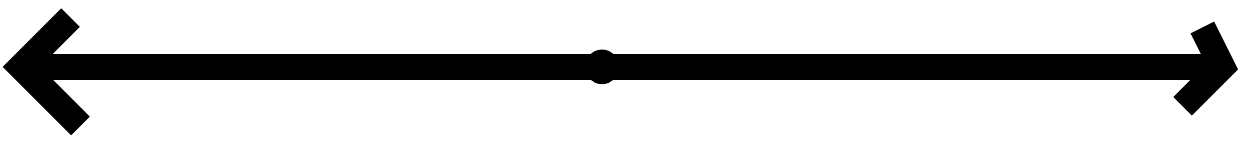}
\end{center}
\end{figure}

That is, if we want to use polynomials, we may   need to prove a more difficult conjecture,
in the sense that we will need to prove non-membership in a larger set.

Given a subset $Z$ of a vector space $U$,   the
{\it ideal of $Z$}, \index{ideal}  denoted $I(Z)$, is the set of all polynomials vanishing
at all points of $Z$. The 
{\it Zariski closure}\index{Zariski closure}  of $Z$, denoted $\ol{Z}$,  is  the set
of $u\in U$ such that $P(u)=0$ for all   $P\in I(Z)$.
The common zero set of a collection of polynomials (such as $\ol{Z}$) is called an
{\it algebraic variety}.

\begin{conjecture} [Mulmuley-Sohoni  \cite{MS1}]   Let $n(m)$ be a polynomial of $m$. Then there exists an $m_0$ such that
for all $m>m_0$,  $\ell^{n-m}\tperm_m \not\in \ol{\tend(\BC^{n^2})\cdot \tdet_n}$.
\end{conjecture}

How serious a problem is the issue of Zariski closure? Does it really change Valiant's conjecture?
  
\smallskip
Mulmuley conjectures  \cite{MS5} that  indeed it does. Namely,
he conjectures that there are sequences in   the closure of the
sequences of spaces $\tend(\BC^{n^2})\cdot \tdet_n$ that   are {\it not} in $\vp$.
 
\begin{example}
Let $P=x_1^3+x_2^2x_3+x_2x_4^2\in S^3\BC^4$.
Then     $P\not\in \tend(\BC^{n^2})\cdot \tdet_n$ for $n<5$ \cite{2015arXiv150502205A}. 
However, $\ol{\tend(\BC^9) \cdot \tdet_3}\supset S^3\BC^4$, i.e., {\it every} homogeneous polynomial of degree $3$ in $4$ variables is in the
Zariski closure of degenerations of $\tdet_3$, see e.g., \cite{MR1786479}.
\end{example}

 However, Mulmuley also conjectures   \cite{MS5}  that any path to resolving Valiant's conjecture will have to address
    \lq\lq wild\rq\rq\ sequences in the closure, so that the stronger conjecture is the more
 natural one.
Moreover Grochow makes the case \cite{Joshunify} that essentially all lower bounds
in algebraic complexity theory have come from algebraic geometry.

\subsection{Advantage of the stronger conjecture: representation theory}
{\it Representation theory} is the systematic study of symmetry in linear algebra.

The variety $\ol{\tend(\BC^{n^2})\cdot \tdet_n}$ may be realized as an {\it orbit closure}  as follows:
Let $GL_N$ denote the group of invertible $N\times N$ matrices. It acts on   the space of polynomials  $S^d\BC^N$ by  \eqref{gactonp}.  
Any element of $\tend(\BC^N)$ my be described as a limit of elements of
$GL_N$, so the Euclidean closure of $\tend(\BC^N)\cdot \tdet_n$
equals the Euclidean closure of $GL_N\cdot \tdet_n$. 
In general Euclidean and Zariski closure can be quite
different (e.g. the Zariski closure of $\{ (z,w)\mid z=0,\ w\in \BZ\}\subset \BC^2$
is the line but this set is already Euclidean closed). However, in this situation Euclidean closure equals Zariski closure (see \cite[Thm. 2.33]{MR1344216}),
so we have the following equality of Zariski closures:
$$\ol{GL_{n^2}\cdot \tdet_n}=\ol{\tend(\BC^{n^2})\cdot \tdet_n}.
$$

  Substantial techniques have been developed to study orbits and their closures.

Let 
$$\hDet_n:=\ol{GL_{n^2}\cdot  \tdet_n }
$$ 
and let 
$$\hPerm^m_n:=\ol{GL_{n^2}\cdot  \ell^{n-m}\tperm_m }.
$$


\begin{conjecture}\label{msmainconj}  \cite{MS1}  Let $n=m^c$ for any constant
$c$. Then for all sufficiently large $m$, 
$$
\hPerm^m_n\not\subset \hDet_n.
$$
\end{conjecture}

 Conjecture \ref{msmainconj} would imply Conjecture \ref{valperdet}.
In \S\ref{MSprogramsect},  I describe the program to use representation theory to prove Conjecture \ref{msmainconj}.

\section{State of the art for conjecture \ref{msmainconj}: classical algebraic geometry}\label{LMRsect}

\subsection{Classical algebraic geometry detour: B. Segre's dimension formula}\label{katzformula}
In algebraic geometry  it is more convenient to work in projective space. 
(From a complexity perspective it is also natural, as changing a function by
a scalar will not change its complexity.) If $W$ is a vector space then
$\BP W$ is the associated projective space of lines through the origin:
$\BP W=(W\backslash 0)/\sim$ where $w_1\sim w_2$ if $w_1=\l w_2$ for some 
nonzero complex number $\l $. Write $[w]\in \BP W$ for
the equivalence class of $w\in W\backslash  0$ and if $X\subset \BP W$, let $\hat X\subset W$ denote
the corresponding cone in $W$. Define  $\ol{X}=\pi(\ol{\hat X})$, the Zariski closure
of $X$.

If $X\subset \BP W$ is a hypersurface,
let $X\dual \subset \BP W^*$ denote its Gauss image, which is called its {\it dual variety}.  If $X$ is an irreducible algebraic
variety, $X\dual$ will be too. More precisely, $X\dual$ is  the Zariski
closure of the set of {\it conormal lines} to smooth points of $X$. Here, if 
$\hat T_xX\subset W$ denotes the tangent space to the cone over $X$,
the conormal space is
$N^*_xX=(\hat T_xX)\upperp \subset W^*$.

\begin{prop}[B. Segre \cite{MR0041481}] Let $P\in S^dW^*$ be irreducible
and let $d\geq 2$.   Then for a Zariski open subset of
points $[x]\in   \tzeros(P)$,
$$\dim \tzeros(P)\dual = \trank (Hess(P)(x^{d-2}))-2.$$
\end{prop}

Here $(Hess(P)(x^{d-2}))\in S^2W^*$ is the Hessian matrix of second partial derivatives of $P$ evaluated
at $x$. Note that the right hand side involves second derivative information, and the
left hand side involves the dual variety (which is first derivative information from
$\tzeros(P)$), and its dimension, which is a first derivative computation on the dual
variety, and therefore a second derivative computation on $\tzeros(P)$.

\begin{proof}  For a homogeneous polynomial $P\in S^dW^*$, write $\ol{P}$ when we consider $P$ as a $d$-multi-linear form.
Let \linebreak
$x\in \hzeros(P)\subset W$ be a smooth  point, so $P(x)=\ol{P}(x\hd x)=0$ and
$dP_x=\ol{P}(x\hd x, \cdot)\neq 0$. Take $h=dP_x\in W^*$, so
$[h]\in \tzeros(P)\dual$.
Consider a curve $h_t\subset \hzeros(P)\dual $ with $h_0=h$. There must be a corresponding
  curve $x_t\in \hzeros(P)$ such that
$h_t=\ol{P}(x_t\hd x_t, \cdot)$ and thus its derivative
is $h_0'=\ol{P}(x^{d-2},x_0' , \cdot)$. 
The dimension of $\hat T_{h}\tzeros(P)\dual$ is then the rank of $Hess(P)(x^{d-2})=\ol{P}(
x^{d-2}, \cdot, \cdot)$ minus one
(we subtract one because $x_0'=x$ is in the
kernel of $Hess(P)(x^{d-2})$). Finally  $\tdim X=\tdim \hat T_xX-1$.
\end{proof}

\subsection{First steps towards equations}\label{firstesteps}
Segre's formula implies, for $P\in S^dW^*$, that   $\tdim \tzeros(P)\dual\leq k$ if and only if,
for all $w\in W$, letting $G(q,W)$ denote the Grassmannian of $q$-planes through the origin in $W$, 
$$
P(w)=0\  \  \Rrightarrow \  \  
 \tdet_{k+3} (Hess(P)(w^{d-2})|_F) =0 \ \forall F\in G(k+3,W). 
$$
  
Equivalently (assuming $P$ is  irreducible), for any   $F\in G(k+3,W)$, the polynomial $P$ must divide \linebreak
$\det_{k+3} (Hess(P) |_F)
\in S^{(k+3)(d-2)}W^*$.

Thus to find polynomials on $S^dW^*$ characterizing hypersurfaces with degenerate duals, we need polynomials
that detect if a polynomial $P\in S^dW^*$ divides a polynomial $Q\in S^eW^*$. 
Now $P$ divides $Q$ if and only if $Q\in P\cdot S^{e-d}W^*$, i.e.,
letting $x^{I_j}$ be a basis of $S^{e-d}W^*$ and let $\ww$ denote  exterior (wedge)
product, 
\be\label{theeqns}
x^{I_1}P\ww\cdots\ww x^{I_D}P\ww Q=0.
\ene
 Let $\tdim W=N$ and let  $\Dual_{k,d,N}\subset \BP S^dW^*$ denote the zero set of the  equations \eqref{theeqns} in the coefficients of $P$ taking $Q=\det_{k+3} (Hess(P) |_F)$.
By our previous discussion $[\tdet_n]\in \Dual_{2n-2,n,n^2}$. 

\subsection{The lower bound on $\ol{dc}(\tperm_m)$}\label{lowerdcp}
When
\be\label{goodpermpt}
x=\begin{pmatrix}1- m & 1 &\cdots &1\\
1  & 1 &\cdots &1\\
\vdots  & \vdots  &\cdots &\vdots \\
1  & 1 &\cdots &1
\end{pmatrix},
\ene
a short calculation shows that     $Hess(\tperm_{m})(x^{m-2})$  is of maximal rank.  This   fills in the missing step of the proof of Theorem \ref{mrthm}. 
 Moreover,    
   if one works over $\BR$, then the Hessian has a signature. For
$\tdet_n$, this signature is $(n-1,n-1)$, but for the permanent the signature
on an open subset is   at least $(m^2-2m+1,2m-3)$, thus:

\begin{theorem}[Yabe \cite{DBLP:journals/corr/Yabe15}] $dc_{\BR}(\tperm_m)\geq m^2-2m+2$.
\end{theorem}

Were we to just consider $\tperm_m$ as a polynomial in more variables, the rank of the Hessian
would not change. However, we are also adding padding, which could {\it a priori} have
a negative effect on the rank of the Hessian. Fortunately, as was shown in \cite{MR3048194} it does not, and we conclude:

\begin{theorem}\cite{MR3048194}\label{dcbarbnd} $\Perm^m_n\not\subset \Dual_{2n-2,n,n^2}$ when $m<\frac {n^2}2$. 
In particular, when $m<\frac {n^2}2$, 
$\Perm^m_n\not\subset \Det_n$.
\end{theorem}
 
On the other hand, since  cones have degenerate duals,  
$\ell^{n-m}\tperm_m \in \Dual_{2n-2,n,n^2}$ whenever $m\geq \frac{n^2}2$.

 In \cite{MR3048194} it was also shown that $\Dual_{k,d,N}$ intersected
 with the set of irreducible hypersurfaces is exactly the set (in $\BP S^dW^*$) of irreducible
 hypersurfaces of degree $d$ in $\BP W$ with dual varieties
 of dimension $k$, which solved a classical question in algebraic geometry.

\section{Necessary conditions for modules of polynomials to be useful
for GCT}\label{gctusefulsect}

Fixing a linear inclusion $\BC^{m^2+1}\subset \BC^{n^2}$, the polynomial $\ell^{n-m}\tperm_m\in S^n\BC^{n^2}$ has evident
pathologies: it is {\it padded}, that is divisible by a large power of a linear form, and its zero set is
a {\it cone} with a $n^2-m^2-1$ dimensional vertex, that is, it only uses $m^2+1$ of the $n^2$ variables
in an expression in good coordinates. 
To separate $\ell^{n-m}\tperm_m$ from $\tdet_n$,
one must  look for modules   in $I(\Det_n)$ that do not vanish automatically
on equations of hypersurfaces with these pathologies. It is easy to determine
such modules with representation theory. Before doing so, 
I first  review  the irreducible representations of the general linear group.

\subsection{$GL(V)$-modules}\label{glvmods}
Let $V$ be a complex vector space of dimension $\bv$.
The irreducible representations of $GL(V)$ are indexed by
sequences of integers $\pi=(p_1\hd p_{\bv})$ with 
$p_1\geq \cdots \geq p_{\bv}$ and 
the corresponding module is denoted $S_{\pi}V$. The representations
occurring in the tensor algebra of $V$ are those with $p_{\bv}\geq 0$, i.e.,
by partitions. For a partition $\pi$, let $\ell(\pi)$ denote its {\it length}, the smallest
$s$ such that $p_{s+1}=0$. 
In particular $S_{(d)}V=S^dV$, and $S_{(1\hd 1)}=: S_{(1^d)}V=\La dV\subset V^{\ot d}$, 
the skew-symmetric tensors. 

One
 way to construct $S_{\pi}V$, where
  $\pi=(p_1\hd p_{\bv})$ and its conjugate partition is
$\pi'=(q_1\hd q_{p_1})$, is to form a projection
operator from $V^{\ot |\pi|}$ by first projecting to
$\La{q_1}V\otc \La{q_{p_1}}V$ by skew-symmetrizing and
then re-ordering and projecting the image to
$S^{p_1}V\otc S^{p_{\bv}}V $. In particular if an element of
$V^{\ot |\pi|}$ lies in some $W^{\ot |\pi|}$ for some $W\subset V$  with $\tdim W<q_1$, then
it will map to zero.

\subsection{Polynomials useful for GCT} \label{gctusesubsect}

To be useful for GCT, a module of polynomials should not vanish identically
on  cones or on polynomials that are divisible by a large power of a linear form.
The equations for the variety of polynomials whose
zero sets are   cones are well known -- they are all modules
where the length of the partition is longer than the
number of variables needed to define the polynomial. 

 
 \begin{prop}\label{gctmainthm}\cite{MR3169697} Necessary conditions for a 
   module  $S_{ {\pi }}\BC^{n^2} \subset I_d(\Det_n)$ to  
not vanish identically on polynomials in $m^2$ variables padded by $\ell^{n-m}$
are
\begin{enumerate}
\item $\ell(\pi )\leq m^2+1 $, 
\item If $\pi =(p_1\hd p_{t})$, then $p_1\geq d(n-m)$.
\end{enumerate}
Moreover, if $p_1\geq \tmin\{d(n-1),dn-m\}$, then the necessary   conditions    are also sufficient.
In particular, for $p_1$ sufficiently large, these conditions depend only on the partition
$\pi$, not how the module $S_{\pi}\BC^{n^2}$  is realized as a space of polynomials.
\end{prop}

\section{The  program to find modules in $I[\Det_n]$ via representation theory}\label{MSprogramsect}
In this section I present   the program initiated in \cite{MS1}  
and developed   in \cite{MR2861717,MS2}
to find modules in the ideal of $\Det_n$.

\subsection{Preliminaries} Let $W=\BC^{n^2}$ and consider $\tdet_n\in S^nW^*$. 
Define $\BC[  \Det_n]:=Sym(S^nW )/I(\Det_n)$, the homogeneous coordinate
ring of   $\hDet_n$. This is the space of polynomial functions
on $\hDet_n$ inherited from polynomials on the ambient space $S^nW$.

Since $Sym(S^nW )$ and $ I(\Det_n)$ are $GL(W)$-modules,
so is $\BC[  \Det_n]$, and since $GL(W)$ is reductive  (a complex
algebraic group $G$ is reductive  if $U\subset V$
is a $G$-submodule of a $G$-module $V$, there exists a complementary 
$G$-submodule $U^c$ such that
$V=U\op U^c$) we obtain
the   splitting as a $GL(W)$-module:
$$
Sym(S^nW )= I(\Det_n)\op \BC[  \Det_n]. 
$$
In particular, if a module $S_{\pi}W $ appears in $Sym(S^nW )$
and it does not appear in $\BC[\Det_n]$, it must appear in $I(\Det_n)$.

For those not familiar with the ring of regular functions
on an affine algebraic variety, consider
$GL(W)\subset \BC^{n^2+1}$ as the 
   subvariety of $\BC^{n^2+1}$, with coordinates
$(x^i_j,t)$  given
by the equation $t\tdet(x)=1$, and $\BC[GL(W)]$ can be defined to be the
restriction of polynomial functions on $\BC^{n^2+1}$ to this subvariety.
  Then $\BC[GL(W)\cdot \tdet_n]=\BC[GL(W)/G_{\tdet_n}]$ can
be defined as the subring of $G_{\tdet_n}$-invariant functions $\BC[GL(W)]^{G_{\tdet_n}}$.
Here $G_{\tdet_n}:=\{ g\in GL(W)\mid g\cdot \tdet_n= \tdet_n\}\simeq SL_n\times SL_n \ltimes \BZ_2$.  
 A nice proof of this result (originally due to Frobenius \cite{Frobdet}) is due to Dieudonn\'e \cite{MR0029360} (see \cite{MR3343444}  for an
exposition). It relies on the fact that, in analogy with a smooth quadric hypersurface, there are two families of
maximal linear spaces on the Grassmannian $G(n^2-n,\BC^n\ot \BC^n)$ with prescribed dimensions of their intersections.
One then uses that the group action must preserve these intersection properties.

  There is an injective map
$$
\BC[ \Det_n]\ra \BC[GL(W)\cdot \tdet_n]
$$
given by restriction of functions. The map is an injection because
any function identically zero on a Zariski open subset of an
irreducible  variety
is identically zero on the variety. 
The algebraic Peter-Weyl theorem below gives a description of
the $G$-module structure of 
$\BC[G/H]$ when $G$ is a reductive 
algebraic group and $H$ is a subgroup.

\smallskip

\noindent  {\bf Plan of \cite{MS1,MS2}}:  {\it Find a module $S_{\pi}W $ not appearing
in $\BC[GL(W)/G_{\tdet_n}]$ that does appear in $Sym(S^nW )$.}

\smallskip

By the above discussion such a module must appear in $I(\Det_n)$.

   One might object that the coordinate rings of different orbits could coincide,
or at least be very close. Indeed this is the case for generic polynomials, but
in GCT one generally restricts to polynomials whose symmetry groups are not
only \lq\lq large\rq\rq , but they {\it characterize} the orbit as follows:

\begin{definition}\label{stabdef} Let $V$ be
a $G$-module. A point  $P\in  V$ is {\it characterized by} its stabilizer $G_P$ if any $Q\in  V$ with $G_Q\supseteq G_P$ is of the
form $Q=cP$ for some constant $c$.
\end{definition}

 One can think of polynomial sequences that are complete for their
complexity classes and    are characterized by their stabilizers as \lq\lq best\rq\rq\  representatives of their class.
Corollary \ref{orbitcor}  will imply that if $P\in S^dV$ is characterized
by its stabilizer, the coordinate ring of its $G$-orbit is unique as a module among orbits of points in $V$.


\subsection{The algebraic Peter-Weyl theorem}\label{algpwsect}
Let $G$ be a complex reductive algebraic group (e.g. $G=GL(W)$), and let
$V$ be an irreducible  $G$-module. Given $v\in V$ and $\a\in V^*$, define
a function $f_{v,\a}: G\ra \BC$ by $f_{v,\a}(g)=\a(g\cdot v)$. 
These are regular functions and it is not hard to see   one obtains an inclusion $V\ot V^*\subset \BC[G]$.
Such functions are called {\it matrix coefficients} as if one takes
bases, these functions  are spanned by the elements of the matrix $\rho(g)$,
where $\rho: G\ra GL(V)$ is the representation.
 In fact the matrix coefficients span $\BC[G]$: 

\begin{theorem}\label{algpw}[Algebraic Peter-Weyl theorem]  Let $G$ be a reductive algebraic group. Then there are only countably many non-isomorphic irreducible
finite dimensional  $G$-modules.
Let $\Lambda_G^+$ denote a set indexing the irreducible $G$-modules,  
and let $V_{\l}$ denote  the irreducible module associated to $\l\in \Lambda_G^+$.
    Then, as a $G\times G$-module
$$
\BC[G]=\bigoplus_{\l \in \Lambda_G^+} V_{\l}\ot  V_{\l}^* .
$$
\end{theorem}

For a proof and discussion, see  e.g. \cite{MR2265844}.

\begin{corollary} \label{orbitcor} Let $H\subset G$ be a   closed  subgroup. Then, as a $G$-module, 
$$
\BC[G/H]=\BC[G]^H=\bigoplus_{\l \in \Lambda_G^+} V_{\l}\ot (V_{\l}^*)^H = \bigoplus_{\l \in \Lambda_G^+} V_{\l}^{\oplus \tdim (V_{\l}^*)^H}.
$$
\end{corollary}
Here $G$ acts on the $V_{\l}$ and
$(V_{\l}^*)^H$ is just a vector space whose dimension records the multiplicity of $V_{\l}$ in $\BC[G/H]$.

Corollary \ref{orbitcor}  motivates the study of polynomials characterized by their stabilizers: if $P\in  V$ is characterized by its stabilizer, 
then $G\cdot P$  is the unique
orbit in $ V$ with   coordinate ring  isomorphic to $\BC[G\cdot P]$ as a $G$-module. 
Moreover, for any $Q\in V$ that is not a multiple of $P$, $\BC[G\cdot Q]\not\subset \BC[G\cdot P]$.

\subsection{Schur-Weyl duality} 
The space $V^{\ot d}$ is acted on by $GL(V)$ and $\FS_d$ (permuting the factors),
and these actions commute so we may decompose it as $GL(V)\times \FS_d$-module.
The decomposition is
$$
V^{\ot d}=\bigoplus_{\pi\mid |\pi|=d} S_{\pi}V\ot [\pi]
$$
where $[\pi]$ is the irreducible $\FS_d$-module associated to
the partition $\pi$, see e.g. \cite{MR1354144}. This gives us a second definition of $S_{\pi}V$
when $\pi$ is a partition: $S_{\pi}V=\thom_{\FS_d}([\pi],V^{\ot d})$.

\subsection{The coordinate ring of $GL(W)\cdot \tdet_n$}\label{detoring}

Let   $E,F\simeq \BC^n$. We first compute the $SL(E)\times SL(F)$-invariants
in $S_{\pi}(E\ot F)$ where $|\pi|=d$.
As a  $GL(E)\times GL(F)$-module, since $(E\ot F)^{\ot d}=E^{\ot d}\ot F^{\ot d}$, 
\begin{align*}
S_{\pi}(E\ot F)&=\thom_{\FS_d}([\pi], E^{\ot d}\ot F^{\ot d})\\
&= \thom_{\FS_d}([\pi], (\bigoplus_{|\mu|=d} [\mu]\ot S_{\mu}E) \ot (\bigoplus_{|\nu|=d} [\nu]\ot S_{\nu}F  ))\\
&= \bigoplus_{|\mu|=|\nu|=d }\thom_{\FS_d}([\pi],   [\mu]\ot [\nu]) \ot S_{\mu}E  \ot   S_{\nu}F  
\end{align*}
The vector space $\thom_{\FS_d}([\pi],   [\mu]\ot [\nu])$ simply records the multiplicity of
$S_{\mu}E  \ot   S_{\nu}F $ in $S_{\pi}(E\ot F)$. 
The integers \linebreak
$k_{\pi  \mu \nu}=\tdim \thom_{\FS_d}([\pi],   [\mu]\ot [\nu])$   are called  {\it  Kronecker coefficients}.

Now  $S_{\mu}E$ is a trivial $SL(E)$ module if and only if $\mu=(\d^n)$ for some
$\d \in \BZ$. Thus so far,  we are reduced to studying the Kronecker coefficients $k_{\pi \d^n \d^n}$. Now
  take the $\BZ_2$ action given by exchanging $E$ and $F$ into account. 
Write $[\mu]\ot [\mu]=S^2[\mu]\op \La 2[\mu]$. The first module will be invariant under  $\BZ_2=\FS_2$, and the second
will transform its sign under the transposition.
So define the {\it symmetric Kronecker coefficients} $sk^{\pi}_{\mu \mu}:= \tdim (\thom_{\FS_d}([\pi],   S^2[\mu]))$.  For a $GL(V)$-module $M$, write
$M_{poly}$ for the submodule consisting of isotypic components of
modules $S_{\pi}V$ where $\pi$ is a partition.

We conclude: 

\begin{prop}\label{peterrefx}\cite{MR2861717} Let $W=\BC^{n^2}$. 
The polynomial part of the coordinate ring of the $GL(W)$-orbit  of $\tdet_n\in S^nW $ is
$$
\BC[GL(W)\cdot  \tdet_n]_{poly}=
 \bigoplus_{d\in \BZ_+}\bigoplus_{\pi \, \mid \, |\pi|=n d }(S_{\pi}W^*)^{\op sk^\pi_{d^n d^n}}.
 $$
  \end{prop}

  


\section{Asymptotics of plethysm and Kronecker coefficients via geometry}\label{kronplethsect}

The above discussion can be summarized as:

\smallskip

\noindent {\bf Goal:} Find partitions $\pi$   satisfying  $ \tmult(S_{\pi}W, S^d(S^nW))\neq 0$,
$sk^\pi_{d^n d^n}=0$, 
   have few parts,  and first part large. 
   
   \smallskip

  Kronecker coefficients and the plethysm coefficients \linebreak
  $\tmult(S_{\pi}W, S^d(S^nW))$ have been well-studied in both the geometry and combinatorics
  literature.  
 I briefly discuss  a geometric method of L. Manivel and J. Wahl \cite{MR1132139,MR1465785,MR1651092,2014arXiv1411.3498M}
  based on the {\it Borel-Weil theorem}  that
 realizes modules as spaces of sections of vector bundles on homogeneous varieties.  Advantages of the method
 are: (i) the vector bundles come with filtrations that allow one to organize information,
 (ii)  
 the sections of the associated graded bundles can be computed explicitly, giving one upper bounds for
 the coefficients, and (iii) Serre's theorem on the vanishing of sheaf cohomology tells one that the
 upper bounds are achieved asymptotically.

A basic, if not {\it the} basic problem in representation theory is: given
a group $G$, an irreducible $G$-module $U$, and a subgroup
$H\subset G$, decompose $U$ as an $H$-module. The determination
of Kronecker coefficients
can be phrased this way with $G=GL(V\ot W)$, $U=S_{\l}(V\ot W)$
and $H=GL(V)\times GL(W)$. The determination of plethysm coefficients
may be phrased as  the case $G=GL(S^nV)$, $U=S^d(S^nV)$
and $H=GL(V)$.

 I focus on plethysm coefficients.
We want to decompose $S^d(S^nV)$ as a $GL(V)$-module, or more precisely,
to obtain qualitative asymptotic information about this decomposition.
Note  that $S^{dn}V\subset S^d(S^nV)$ with multiplicity
one. Let $x_1\hd x_{\bv}$ be a basis of $V$, so $((x_1)^n)^d$ is
the highest highest weight vector in $S^d(S^nV)$. (A vector $v\in V$ is a highest
weight vector for $GL(W)$ if $B[v]=[v]$ where $B\subset GL(W)$ is the subgroup
of upper triangular matrices. There is a partial order on the set of highest weights.) 
Say $S_{\pi}V\subset S^d(S^nV)$ is realized with highest weight vector
$$
\sum_I c^I (x_{i_{11}}\cdots x_{i_{1n}})\cdots
(x_{i_{d1}}\cdots x_{i_{dn}})
$$
for some coefficients $c^I$, where $I=\{ i_{s,\a}\}$.
Then
$$
\sum_I c^I (x_1)^n (x_{i_{11}}\cdots x_{i_{1n}})\cdots
(x_{i_{d1}}\cdots x_{i_{dn}})\in S^{d+1}(S^nV)
$$
is a vector of weight $(n)+\pi$, and is   a highest weight vector.
Similarly
$$
\sum_I c^I (x_1x_{i_{11}}\cdots x_{i_{1n}})\cdots
(x_1x_{i_{d1}}\cdots x_{i_{dn}})\in S^{d}(S^{n+1}V)
$$
is a vector of weight $(d)+\pi$, and is   a highest weight vector.
This already shows qualitative behavior if we allow the first part of
a partition to grow:
\begin{prop}\cite{MR1465785}\label{manprop}
Let $\mu$ be a fixed partition. Then  $\tmult(S_{(dn-|\mu|,\mu)}, S^d(S^nV))$ is a non-decreasing  
function of both $d$ and $n$.
\end{prop}

One way to view what we just did was to write
$V=x_1\op T$,  so
\be\label{oscver}S^n(x_1\op T)=\bigoplus_{j=0}^nx_1^{n-j}\ot S^jT.
\ene
Then  decompose  the $d$-th symmetric power of $S^n(x_1\op T)$  and examine
the stable behaviour as we increase $d$ and $n$. One could think of
the decomposition \eqref{oscver}  as the osculating sequence of
the $n$-th Veronese embedding of $\BP V$ at 
$[x_1^n]$ and the further decomposition as the osculating sequence of the $d$-th Veronese
re-embedding of the ambient space   refined by \eqref{oscver}.

For Kronecker coefficients and more general decomposition problems
the situation is more complicated in that the ambient space is
no longer be projective space,  but a homogeneous variety, and instead
of an osculating sequence, one   examines jets of sections of a vector bundle. 
As mentioned above, in this situation one gets the bonus of vanishing theorems.
For example, with the use of vector bundles, Proposition \ref{manprop} can be strengthened to
say that the multiplicity is eventually constant and state for which    $d,n$  
this constant multiplicity is achieved.

\subsection*{Acknowledgements} I thank Jesko H\"uttenhain for drawing the   pictures of surfaces, and H. Boas and J. Grochow for extensive suggestions for improving the exposition.

 \printnotes

\bibliographystyle{amsalpha}
\bibliography{Lmatrix}

\def\cdprime{$''$} \def\cprime{$'$} \def\cprime{$'$} \def\cprime{$'$}
  \def\Dbar{\leavevmode\lower.6ex\hbox to 0pt{\hskip-.23ex \accent"16\hss}D}
  \def\cprime{$'$} \def\cprime{$'$} \def\cdprime{$''$} \def\cprime{$'$}
  \def\cprime{$'$} \def\Dbar{\leavevmode\lower.6ex\hbox to 0pt{\hskip-.23ex
  \accent"16\hss}D} \def\cprime{$'$} \def\cprime{$'$} \def\cprime{$'$}
  \def\cprime{$'$} \def\Dbar{\leavevmode\lower.6ex\hbox to 0pt{\hskip-.23ex
  \accent"16\hss}D} \def\cprime{$'$} \def\cprime{$'$}
\providecommand{\bysame}{\leavevmode\hbox to3em{\hrulefill}\thinspace}
\providecommand{\MR}{\relax\ifhmode\unskip\space\fi MR }
\providecommand{\MRhref}[2]{%
  \href{http://www.ams.org/mathscinet-getitem?mr=#1}{#2}
}
\providecommand{\href}[2]{#2}
\begin{thebibliography}{KLPSMN09}

\bibitem[ABV15]{2015arXiv150502205A}
J.~{Alper}, T.~{Bogart}, and M.~{Velasco}, \emph{{A lower bound for the
  determinantal complexity of a hypersurface}}, ArXiv e-prints (2015).

\bibitem[Alu]{aluffideg}
Paolo Aluffi, \emph{Degrees of projections of rank loci}, preprint
  arXiv:1408.1702.

\bibitem[BCS97]{BCS}
Peter B{\"u}rgisser, Michael Clausen, and M.~Amin Shokrollahi, \emph{Algebraic
  complexity theory}, Grundlehren der Mathematischen Wissenschaften
  [Fundamental Principles of Mathematical Sciences], vol. 315, Springer-Verlag,
  Berlin, 1997, With the collaboration of Thomas Lickteig. \MR{99c:68002}

\bibitem[Bea00]{MR1786479}
Arnaud Beauville, \emph{Determinantal hypersurfaces}, Michigan Math. J.
  \textbf{48} (2000), 39--64, Dedicated to William Fulton on the occasion of
  his 60th birthday. \MR{1786479 (2002b:14060)}

\bibitem[BLMW11]{MR2861717}
Peter B{\"u}rgisser, J.~M. Landsberg, Laurent Manivel, and Jerzy Weyman,
  \emph{An overview of mathematical issues arising in the geometric complexity
  theory approach to {${\rm VP}\neq{\rm VNP}$}}, SIAM J. Comput. \textbf{40}
  (2011), no.~4, 1179--1209. \MR{2861717}

\bibitem[Cai90]{MR1032157}
Jin-Yi Cai, \emph{A note on the determinant and permanent problem}, Inform. and
  Comput. \textbf{84} (1990), no.~1, 119--127. \MR{MR1032157 (91d:68028)}

\bibitem[Coo71]{cook1971complexity}
Stephen~A Cook, \emph{The complexity of theorem-proving procedures},
  Proceedings of the third annual ACM symposium on Theory of computing, ACM,
  1971, pp.~151--158.

\bibitem[Die49]{MR0029360}
Jean Dieudonn{\'e}, \emph{Sur une g\'en\'eralisation du groupe orthogonal \`a
  quatre variables}, Arch. Math. \textbf{1} (1949), 282--287. \MR{0029360
  (10,586l)}

\bibitem[Fro97]{Frobdet}
G.~Frobenius, \emph{{\"U}ber die {D}arstellung der endlichen {G}ruppen durch
  lineare {S}ubstitutionen}, Sitzungsber Deutsch. Akad. Wiss. Berlin (1897),
  994--1015.

\bibitem[Gat87]{MR922386}
Joachim von~zur Gathen, \emph{Feasible arithmetic computations: {V}aliant's
  hypothesis}, J. Symbolic Comput. \textbf{4} (1987), no.~2, 137--172.
  \MR{MR922386 (89f:68021)}

\bibitem[GHIL]{GHILrigid}
Fulvio Gesmundo, Jonathan Hauenstein, Christian Ikenmeyer, and J.~M. Landsberg,
  \emph{Geometry and matrix rigidity}, to appear in FOCM, arXiv:1310.1362.

\bibitem[Gre14]{Gre11}
Bruno Grenet, \emph{{An Upper Bound for the Permanent versus Determinant
  Problem}}, Theory of Computing (2014), Accepted.

\bibitem[Kar72]{karp1972reducibility}
Richard~M Karp, \emph{Reducibility among combinatorial problems}, Springer,
  1972.

\bibitem[KL14]{MR3169697}
Harlan Kadish and J.~M. Landsberg, \emph{Padded polynomials, their cousins, and
  geometric complexity theory}, Comm. Algebra \textbf{42} (2014), no.~5,
  2171--2180. \MR{3169697}

\bibitem[KLPSMN09]{MR2870721}
Abhinav Kumar, Satyanarayana~V. Lokam, Vijay~M. Patankar, and Jayalal Sarma
  M.~N., \emph{Using elimination theory to construct rigid matrices},
  Foundations of software technology and theoretical computer
  science---{FSTTCS} 2009, LIPIcs. Leibniz Int. Proc. Inform., vol.~4, Schloss
  Dagstuhl. Leibniz-Zent. Inform., Wadern, 2009, pp.~299--310. \MR{2870721}

\bibitem[Lan08]{MR2383305}
J.~M. Landsberg, \emph{Geometry and the complexity of matrix multiplication},
  Bull. Amer. Math. Soc. (N.S.) \textbf{45} (2008), no.~2, 247--284.
  \MR{MR2383305 (2009b:68055)}

\bibitem[Lan15]{MR3343444}
\bysame, \emph{Geometric complexity theory: an introduction for geometers},
  Ann. Univ. Ferrara Sez. VII Sci. Mat. \textbf{61} (2015), no.~1, 65--117.
  \MR{3343444}

\bibitem[LMR13]{MR3048194}
Joseph~M. Landsberg, Laurent Manivel, and Nicolas Ressayre, \emph{Hypersurfaces
  with degenerate duals and the geometric complexity theory program}, Comment.
  Math. Helv. \textbf{88} (2013), no.~2, 469--484. \MR{3048194}

\bibitem[LR15]{LRpermdet}
J.M. Landsberg and Nicolas Ressayre, \emph{Permanent v. determinant: an
  exponential lower bound assuming symmetry and a potential path towards
  valiant's conjecture}, preprint (2015).

\bibitem[Mac95]{MR1354144}
I.~G. Macdonald, \emph{Symmetric functions and {H}all polynomials}, second ed.,
  Oxford Mathematical Monographs, The Clarendon Press Oxford University Press,
  New York, 1995, With contributions by A. Zelevinsky, Oxford Science
  Publications. \MR{1354144 (96h:05207)}

\bibitem[Man97]{MR1465785}
Laurent Manivel, \emph{Applications de {G}auss et pl\'ethysme}, Ann. Inst.
  Fourier (Grenoble) \textbf{47} (1997), no.~3, 715--773. \MR{MR1465785
  (98h:20078)}

\bibitem[Man98]{MR1651092}
\bysame, \emph{Gaussian maps and plethysm}, Algebraic geometry ({C}atania,
  1993/{B}arcelona, 1994), Lecture Notes in Pure and Appl. Math., vol. 200,
  Dekker, New York, 1998, pp.~91--117. \MR{MR1651092 (99h:20070)}

\bibitem[{Man}14]{2014arXiv1411.3498M}
L.~{Manivel}, \emph{{On the asymptotics of Kronecker coefficients}}, ArXiv
  e-prints (2014).

\bibitem[MM61]{MR0147488}
Marvin Marcus and Henryk Minc, \emph{On the relation between the determinant
  and the permanent}, Illinois J. Math. \textbf{5} (1961), 376--381.
  \MR{0147488 (26 \#5004)}

\bibitem[MN]{MS5}
Ketan~D. Mulmuley and H.~Narayaran, \emph{Geometric complexity theory {V}: On
  deciding nonvanishing of a generalized {L}ittlewood-{R}ichardson
  coefficient}, Technical Report TR-2007-05, computer science department, The
  University of Chicago, May, 2007.

\bibitem[MR04]{MR2126826}
Thierry Mignon and Nicolas Ressayre, \emph{A quadratic bound for the
  determinant and permanent problem}, Int. Math. Res. Not. (2004), no.~79,
  4241--4253. \MR{MR2126826 (2006b:15015)}

\bibitem[MS01]{MS1}
Ketan~D. Mulmuley and Milind Sohoni, \emph{Geometric complexity theory. {I}.
  {A}n approach to the {P} vs.\ {NP} and related problems}, SIAM J. Comput.
  \textbf{31} (2001), no.~2, 496--526 (electronic). \MR{MR1861288
  (2003a:68047)}

\bibitem[MS08]{MS2}
\bysame, \emph{Geometric complexity theory. {II}. {T}owards explicit
  obstructions for embeddings among class varieties}, SIAM J. Comput.
  \textbf{38} (2008), no.~3, 1175--1206. \MR{MR2421083}

\bibitem[Mum95]{MR1344216}
David Mumford, \emph{Algebraic geometry. {I}}, Classics in Mathematics,
  Springer-Verlag, Berlin, 1995, Complex projective varieties, Reprint of the
  1976 edition. \MR{1344216 (96d:14001)}

\bibitem[Pro07]{MR2265844}
Claudio Procesi, \emph{Lie groups}, Universitext, Springer, New York, 2007, An
  approach through invariants and representations. \MR{MR2265844 (2007j:22016)}

\bibitem[Seg10]{csegre}
C.~Segre, \emph{Preliminari di una teoria delle variet\`a luoghi di spazi},
  Rend. Circ. Mat. Palermo (1910), no.~XXX, 87--121.

\bibitem[Seg51]{MR0041481}
Beniamino Segre, \emph{Bertini forms and {H}essian matrices}, J. London Math.
  Soc. \textbf{26} (1951), 164--176. \MR{0041481 (12,852g)}

\bibitem[Sip92]{Sipser}
Michael Sipser, \emph{The history and status of the p versus np question}, STOC
  '92 Proceedings of the twenty-fourth annual ACM symposium on Theory of
  computing (1992), 603--618.

\bibitem[Spi79]{MR532832}
Michael Spivak, \emph{A comprehensive introduction to differential geometry.
  {V}ol. {III}}, second ed., Publish or Perish Inc., Wilmington, Del., 1979.
  \MR{MR532832 (82g:53003c)}

\bibitem[Val79]{MR526203}
L.~G. Valiant, \emph{The complexity of computing the permanent}, Theoret.
  Comput. Sci. \textbf{8} (1979), no.~2, 189--201. \MR{MR526203 (80f:68054)}

\bibitem[vzG87]{MR910987}
Joachim von~zur Gathen, \emph{Permanent and determinant}, Linear Algebra Appl.
  \textbf{96} (1987), 87--100. \MR{MR910987 (89a:15005)}

\bibitem[Wah91]{MR1132139}
Jonathan Wahl, \emph{Gaussian maps and tensor products of irreducible
  representations}, Manuscripta Math. \textbf{73} (1991), no.~3, 229--259.
  \MR{1132139 (92m:14066a)}

\bibitem[Yab15]{DBLP:journals/corr/Yabe15}
Akihiro Yabe, \emph{Bi-polynomial rank and determinantal complexity}, CoRR
  \textbf{abs/1504.00151} (2015).

\end{thebibliography}

\begin{info}
Joseph (J.M.)   Landsberg [\url{jml@math.tamu.edu}] is a professor of mathematics at Texas A\&M University. He has broad research interests, most recently applying
geometry and representation theory to questions in theoretical computer
science. He is co-author (with T. Ivey)  of {\it Cartan for Beginners} (AMS GSM 61) and author of {\it Tensors: Geometry and Applications}
(AMS GSM 128). In the fall of 2014, Landsberg served as Chancellor's  Professor at
the Simons Institute for the Theory of Computing, UC Berkeley.
\end{info}

\end{document}